\documentclass{article}
\usepackage{amssymb}
\usepackage{amsmath}
\usepackage{amsfonts}

\newtheorem{theorem}{Theorem}

\newtheorem{definition}[theorem]{Definition}

\newtheorem{proposition}[theorem]{Proposition}

\input{tcilatex}
\begin{document}

\title{On a simple proof of the Poincare Conjecture}
\author{Erc\"{u}ment Orta\c{c}gil}
\maketitle

\begin{abstract}
We outline a simple proof of PC without surgeries using the homogeneous flow
introduced in [O].
\end{abstract}

\section{The main idea of the new proof}

The Ricci flow introduced by Hamilton deforms an arbitrary metric in the
direction of its Ricci curvature. The difficulty of attacking the Poincare
Conjecture (PC) using the Ricci flow is the lack of any canonical metric on
a homotopy 3-sphere and the fact that the choice of an arbitrary metric as
initial condition produces singularities. In [O] we introduced a flow on
parallelizable manifolds, called the homogeneous flow HF, which strarts with
an arbitrary parallelism as initial condition and flows in the direction of
the curvature of this parallelism. By the theory developed in [O], a
parallelism with vanishing curvature defines a Lie group structure on the
underlying compact and simply connected manifold. This fact points at the
possibility of proving PC by putting a Lie group structure on a homotopy
3-sphere. The key fact, as we will outline in the next section, is the
existence of \textit{canonical parallelisms on a homotopy 3-sphere }$M$
studied first by [H], [A] and [KM] which coincide with the canonical
left-right Lie orbits if $M=S^{3}.$ The main observation of this short
announcement is that HF converges for these canonical initial conditions
whereby proving PC. We will postpone the technical details of the proof to
some future work and content here by outlining some details of the overall
picture.

\section{The framework of HF}

We start by recalling some facts from [O] and refer to [O] for further
details. Let $M$ be a smooth manifold with $\dim M=n$ and $%
F(M)\longrightarrow M$ be the principal \textit{coframe} bundle, i.e., the 
\textit{left} $GL(n,\mathbb{R})$-principal bundle whose fiber over $x\in M$
is the set of all $1$-jets of local diffeomorphisms (which we call $1$%
-arrows) with source at $x\in M$ and target at the origin $o\in \mathbb{R}%
^{n}.$ Therefore, $GL(n,\mathbb{R)}$ is identified with the group of all $1$%
-arrows with source and target at $o$ and acts on $1$-arrows by composition
at the target. Henceforth we always assume that $F(M)\longrightarrow M$
admits a global crossection $w,$ i.e., $M$ is parallelizable. In terms of
coordinates, $w=(w_{j}^{(i)}(x))$ on $(U,x)$ with $\det w\neq 0$ and a
coordinate change $(x)\longrightarrow (y)$ transforms the components of $w$
according to

\begin{equation}
w_{a}^{(i)}(x)\frac{\partial x^{a}}{\partial y^{j}}=w_{j}^{(i)}(y)\text{ \ \
\ \ \ \ }1\leq i\leq n
\end{equation}

Therefore $w$ defines $n$ independent global $1$-forms $w^{(i)}$ on $M,$ $%
1\leq i\leq n,$ which are of the form $w^{(i)}(x)=w_{a}^{(i)}(x)dx^{a}$ on $%
(U,x).$ We called $w$ the structure object in [O] whereas trivialization or
framing are the standard topological terms. Note that $\overline{w}\overset{%
def}{=}w^{-1}$ defines a crossection of the principal \textit{frame }bundle
(a \textit{right }$GL(n,\mathbb{R})$-principal bundle) and we may work
either with $w$ or $\overline{w}.$

Now let $\mathcal{A}\longrightarrow M$ be the group bundle whose fiber $%
\mathcal{A}^{x}$ over $x\in M$ is the group of all $1$-arrows with source
and target at $x\in M.$ A choice of coordinates around $x$ identifies $%
\mathcal{A}^{x}$ with $GL(n,\mathbb{R}).$ Global sections $\Gamma \mathcal{A}
$ of $\mathcal{A}\longrightarrow M$ are called gauge transformations and
they form a group called the gauge group under fiberwise composition. Let $%
\mathcal{W}=\mathcal{W}(M)$ be the set of all structure objects on $M.$ Now $%
\Gamma \mathcal{A}$ acts on $\mathcal{W}$ on the \textit{right }by
composition of $1$-arrows at the source. In coordinates, if $%
w=(w_{j}^{(i)}(x))$ and $a=(a_{j}^{i}(x)),$ then $\left( w\circ a\right)
_{j}^{(i)}(x)\overset{def}{=}w_{s}^{(i)}(x)a_{j}^{s}(x).$ This action is
clearly free. It is also transitive because if $w,z$ are two structure
objects, then $w^{-1}\circ z=a$ is a gauge transformation which maps $w$ to $%
z.$ Therefore this action is simply transitive.

A gauge transformation $a=(a_{j}^{i}(x))$ is positive if $\det $\bigskip $%
(a_{j}^{i}(x))$ is positive for all $x\in M.$ Since a coordinate change
transforms $a_{j}^{i}(x)$ according to $\frac{\partial y^{i}}{\partial x^{t}}%
a_{s}^{t}(x)\frac{\partial x^{s}}{\partial y^{j}}=a_{j}^{i}(y),$ this
definition is independent of coordinates. Since $M$ is parallelizable, it is
clearly orientable. Suppose we fix once and for all an atlas whose
transition functions have positive determinants. A structure object $%
w=(w_{j}^{(i)}(x))$ is positive (with respect to this orientation) if $\det
(w_{j}^{(i)}(x))$ is positive and this definition is independent of the
oriented coordinates by (1). Note that $\lambda w,$ $0\neq \lambda \in 
\mathbb{R}$ is another structure object and $w$ and $-w$ belong to different
orientations if $\dim M$ is odd. Henceforth we will assume that $\mathcal{W}$
and $\Gamma \mathcal{A}$ denote positive structure objects and positive
gauge transformations so that $\Gamma \mathcal{A}$ acts simply transitively
on $\mathcal{W}$ on the right.

Let $GL^{+}(n,\mathbb{R})$ be the group of all $1$-arrows with source and
target at the origin $o\in \mathbb{R}^{n}$ with positive determinant and we
recall that $\mathcal{A}^{x}$ denotes the fiber of $\mathcal{A}%
\longrightarrow M$ over $x\in M$ consisting of all $1$-arrows with source
and target at $x$ with positive determinant. Let $C(M,GL^{+}(n,\mathbb{R}))$
denote the set of smooth functions on $M$ with values in $GL^{+}(n,\mathbb{R}%
).$ Now $C(M,GL^{+}(n,\mathbb{R}))$ acts on $\mathcal{W}$ by composition of $%
1$-arrows on the left as $\left( f\circ w\right) _{j}^{(i)}(x)\overset{def}{=%
}f_{(s)}^{(i)}(x)w_{j}^{(s)}(x)$ and clearly this action is also simply
transitive. Therefore, some $w\in \mathcal{W}$ determines bijections between 
$\mathcal{W},$ $\Gamma \mathcal{A}$ and $C(M,GL^{+}(n,\mathbb{R})).$ More
explicitly, if $w\in \mathcal{W},$ then $w(x)$ defines an isomorphism $%
\overline{w}(x)$

\begin{eqnarray}
\overline{w}(x) &:&\mathcal{A}^{x}\longrightarrow GL^{+}(n,\mathbb{R}) \\
&:&a_{x}\longrightarrow w(x)\circ a_{x}\circ w^{-1}(x)  \notag
\end{eqnarray}

Now we fix $w,$ choose $w^{\prime }\in \mathcal{W}$ arbitrarily and let $%
a\in \Gamma \mathcal{A}$ be the unique gauge transformation satisfying $%
w\circ a=w^{\prime }.$ Now (2) gives the map

\begin{eqnarray}
\overline{w} &:&\mathcal{W\longrightarrow }C(M,GL^{+}(n,\mathbb{R})) \\
&:&w^{\prime }\longrightarrow \{x\longrightarrow w(x)\circ a_{x}\circ
w^{-1}(x)=w^{\prime }(x)\circ w^{-1}(x)\}  \notag
\end{eqnarray}%
which is easily seen to be a bijection. Note that $\overline{w}$ maps $w$ to
the constant map $1\in GL^{+}(n,\mathbb{R}).$ Clearly $C(M,GL^{+}(n,\mathbb{R%
}))$ is a group with the obvious multiplication. We define $w^{\prime }\cdot
w^{\prime \prime }\overset{def}{=}w\circ ab$ where $w^{\prime }=w\circ a$
and $w^{\prime \prime }=w\circ b.$ With this operation $\mathcal{W}$ is
easily seen to be a group with identity $w$ and $\overline{w}$ defined by
(3) becomes an isomorphism. In short, some $w\in \mathcal{W}$ identifies the
groups $\mathcal{W},$ $\Gamma \mathcal{A}$ and $C(M,GL^{+}(n,\mathbb{R})).$

Now let $[M,GL^{+}(n,\mathbb{R})]$ denote the set of (smooth) homotopy
classes of maps. The multiplication on $C(M,GL^{+}(n,\mathbb{R}))$ descends
to $[M,GL^{+}(n,\mathbb{R})]$ turning it into a group. Our purpose is to
define an equivalence relation $\thicksim $ on $\mathcal{W},$ turn $\mathcal{%
W}/\thicksim $ into a group in such a way that the equivalence class $[w]$
of $w$ becomes the identity and defines an isomorphism

\begin{equation}
\overline{\lbrack w]}:\mathcal{W}/\sim \text{ }\longrightarrow \text{ }%
[M,GL^{+}(n,\mathbb{R})]
\end{equation}

For this purpose, let $a(t),$ $0\leq t\leq 1,$ be a smooth $1$-parameter
family of gauge transformation starting from the identity, i.e., $a(0)=Id.$
In coordinates, $a(t)=(a_{j}^{i}(t,x)),$ $a_{j}^{i}(0,x)=\delta _{j}^{i}$
for all $x\in M.$ We called $a(t)$ a deformation in [O]. For $w^{\prime
},w^{\prime \prime }\in \mathcal{W},$ $w^{\prime }\thicksim w^{\prime \prime
}$ if there exists a deformation (or homotopy) $a(t)$ satisfying $w^{\prime
}\circ a(1)=w^{\prime \prime }.$ In other words, $w^{\prime }\circ a(t)$ is
a smooth family of structure objects (= trivializations = framings) starting
from $w^{\prime }$ for $t=0$ and ending at $w^{\prime \prime }$ for $t=1.$
Note that given some $w(t)$ with $w(0)=w_{0},$ if we define

\begin{equation}
w(0)^{-1}\circ w(t)\overset{def}{=}a(t)
\end{equation}%
then $a(t)$ is a gauge transformation for all $t,$ $a(0)=Id$ and

\begin{equation}
w_{0}\circ a(t)=w(t)
\end{equation}

It is easy to check that the multiplication on $\mathcal{W}$ defined above
descends to $\mathcal{W}/\thicksim $ turning it into a group with the
identity $[w]$ and $\overline{w}$ descends to the isomorphism $\overline{[w]}
$ in (4). Note that $\mathcal{W}/\sim $ may be interpreted as the set $\pi
_{0}(\Gamma \mathcal{A})$ of path components of the group $\Gamma \mathcal{A}%
.$ We also observe the crucial fact that the choice of $[w]=\alpha \in $ $%
\mathcal{W}/\sim $ is arbitrary and \textit{there is no canonical choice in
general}. However, canonical orbits \textit{do sometimes exist} which we now
come to.

We recall that from the standpoint of [O], a Lie group structure on $M,$
interpreted as a simply transitive transformation group of $M,$ is a special 
$w.$ In more detail, some $w\in \mathcal{W}$ defines a subgroupoid $%
\varepsilon (M\times M)\subset \mathcal{U}_{1}=$ the universal groupoid of
all $1$-arrows on $M.$ The unique $1$-arrow of $\varepsilon (M\times M)$
from $p$ to $q$ is $w(q)^{-1}\circ w(p).$ If the linear curvature $\mathfrak{%
R}$ (or equivalently the nonlinear curvature $\mathcal{R})$ of $\varepsilon
(M\times M)$ vanishes identically, then the $1$-arrows of $\varepsilon
(M\times M)$ integrate uniquely to a transitive pseudogroup $\mathcal{G}$ on 
$M.$ In this case we call $(M,w)$ a local Lie group (LLG) and denote it by $%
(M,w,\mathcal{G)}.$ Note that this a global concept. If the local
diffeomorphism of $\mathcal{G}$ extend (necessarily uniquely) to global
diffeomorphisms of $M$ (which is the case if $M$ is compact and simply
connected), then $\mathcal{G}$ defines a simply transitive transformation
group of $M.$ Now we can define an abstract Lie group structure on $M$
making $\mathcal{G}$ left or right translations depending on our choice. In
short, a Lie group structure on $M$ is a $w$ with vanishing curvature which
is further globalizable. We refer to Part 1 of [O] for the technical details
of this construction. In this case, we will call $w$ Lie and its equivalence
class $[w]$ a Lie orbit. If $(M_{1},w_{1},\mathcal{G}_{1})$ and $%
(M_{2},w_{2},\mathcal{G}_{2})$ are two Lie groups, then they are isomorphic
if there exists a diffeomorphism $f:M_{1}\longrightarrow M_{2}$ satisfying $%
f\circ \mathcal{G}_{1}\circ f^{-1}=\mathcal{G}_{2}.$ If we denote the Lie
arbit $[w]$ by $\alpha _{1},$ then there is another Lie orbit $\alpha
_{2}\in \mathcal{W}/\sim $ and $\alpha _{1},\alpha _{2}$ represent
left/right or right/left Lie orbits depending on our choice. Note that $%
\alpha _{1}$ and $\alpha _{2}$ need not be the same (see (15) below for $%
M=S^{3}).$

Now let $G$ be an abstract Lie group with the underlying manifold $M$ and
suppose we identify $G$ with the transformation group $LG=\{L_{g}\mid g\in
G\}$ of left translations. The $1$-arrows of the transformations of $LG$
defines a subgroupoid $\varepsilon (M\times M)\subset \mathcal{U}_{1}.$ We
fix some $p\in M$ and consider all $1$-arrows of $\varepsilon (M\times M)$
with source at $x\in M$ and target at $p.$ Fixing a coordinate system around 
$p$ and identifying $p$ with $o\in \mathbb{R}^{n},$ we get a structure
object $w.$ Changing $p$ or the coordinate system around $p$ amounts to
acting with some $g\in GL^{+}(n,\mathbb{R})$ on the $1$-arrows of $w$ at the
target at $o\in \mathbb{R}^{n},$ i.e., $w(x)\longrightarrow g\circ w(x),$ $%
x\in M.$ Note that $w$ and $g\circ w$ define the same groupoid $\varepsilon
(M\times M)$ because $w(q)^{-1}\circ w(p)=\left( g\circ w(p)\right)
^{-1}\circ $ $\left( g\circ w(p)\right) .$ However, acting with $g$ at the
target is the same as acting with $w(x)^{-1}\circ g\circ w(x)$ at the source
since $g\circ w(x)=w(x)\circ \left( w(x)^{-1}\circ g\circ w(x)\right) .$ Now
connecting $g$ with a path to $1\in GL^{+}(n,\mathbb{R})$ gives a
deformation of $w,$ proving that $[w]$ is independent of these choices. Note
that this deformation is very special: it is $\varepsilon $-invariant and
not all deformations are that simple. The same construction works if we
replace $LG$ with $RG.$ To summarize, an abstract Lie group $G$ with the
base manifold $M$ defines \textit{two} \textit{canonical Lie orbits} $\alpha
_{1},\alpha _{2}\in \mathcal{W}/\thicksim $ corresponding \textit{%
canonically }to left and right translations respectively.

Now let $SO(n)\subset $ $GL^{+}(n,\mathbb{R})$ be the subgroup of all $1$%
-arrows with source and target at $o\in \mathbb{R}^{n}$ fixing the canonical
metric $(\delta _{ij})$ and having positive determinant. Then $SO(n)\subset $
$GL^{+}(n,\mathbb{R})$ is a maximal compact subgroup and therefore $GL^{+}(n,%
\mathbb{R})/SO(n)$ is contractible. It follows that $SO(n)$ has the same
homotopy type as $GL^{+}(n,\mathbb{R})$ and therefore the set inclusion $%
C(M,SO(n))\subset C(M,GL^{+}(n,\mathbb{R}))$ induces an isomorphism

\begin{equation}
\lbrack M,GL^{+}(n,\mathbb{R})]\cong \lbrack M,SO(n)]
\end{equation}

Now we recall the homogeneous flow (HF) defined in [O]. Given some $w,$ we
can define the tensors $\mathcal{H}^{(k)}=(\mathcal{H}_{j}^{(k)}),$ $1\leq
k\leq n$ using the linear curvature $\mathfrak{R}=(\mathfrak{R}_{jk,l}^{i}).$
This is analogous to defining the Ricci curvature of the Riemann curvature
tensor. Given some $w_{0}\in \mathcal{W},$ let $w(t)=(w_{j}^{(i)}(t,x))$ be
a flow of structure objects satisfying

\begin{equation}
\frac{\partial w_{j}^{(i)}(t,x)}{\partial t}=\mathcal{H}_{j}^{(i)}(w(t,x))%
\text{ \ \ \ \ \ }w(0)=w_{0}
\end{equation}

If $M$ is compact, then for any initial condition $w_{0}\in \mathcal{W},$ HF
is defined by (8) and has unique solution for $0\leq t\lneqq \epsilon $ for
some $\epsilon .$ Now multiplying both sides of (8) by $\overline{w}%
_{(i)}^{k}(0)$ where $\overline{w}=w^{-1},$ summing over $i$ and using (5),
(6), we rewrite (8) in the form

\begin{equation}
\frac{\partial a_{j}^{k}(t,x)}{\partial t}=\mathcal{H}_{j}^{k}(a(t,x),w_{0})%
\text{ \ \ \ \ }a(0)=Id
\end{equation}

If $a(t)$ is the solution of (9) for $0\leq t\leq t^{\prime }$ for some $%
w_{0},$ note that $w_{0}$ is homotopic to $w_{0}\circ a(t),$ i.e., $%
[w_{0}]=[w_{0}\circ a(t)],$ for all $t\in \lbrack 0,t^{\prime }].$
Therefore, if we plug some inital condition into HF and run HF, we remain
always in the orbit of the initial condition in $\mathcal{W}/\thicksim .$

\begin{definition}
Some $w_{0}\in \mathcal{W}$ is convergent if

i) The unique solution $a(t)$ of (9) is defined for all $t\geq 0.$

ii) There exists a $t^{\prime }\geq 0$ such that $a(t)=a(t^{\prime })$ for
all $t\geq t^{\prime }.$
\end{definition}

If $w_{0}$ is convergent, then $w_{0}\circ a(t)\in \mathcal{W},$ $t\geq
t^{\prime },$ is the limit of $w_{0}.$ The $\mathcal{H}$-tensor of this
limit vanishes since the left hand side of (9) vanishes eventually since $%
a(t)$ becomes constant. The extreme example of convergence is the case of a
Lie group structure (or more generally a LLG) $(M,w_{0},\mathcal{G}).$ Since 
$\mathfrak{R}=0$ implies $\mathcal{H}=0,$ with $w_{0}$ as the initial
condition, $a(t)=Id,$ $t\geq 0$ is the unique solution of (9) and therefore
we may choose $t^{\prime }=0$ in Definition 1.

We now specialize to $M=S^{3}.$ Since $S^{3}$ is compact and $\dim
S^{3}=\dim SO(3)=3,$ the degree of a map $f:S^{3}\longrightarrow SO(3)$ is
defined and furthermore, two such maps are homotopic if and only if they
have the same degree. Therefore we have the isomorphism

\begin{equation}
\deg :[S^{3},SO(3)]=\pi _{3}(SO(3))\cong \mathbb{Z}
\end{equation}

The maps $f:S^{3}\longrightarrow SO(3)$ with $\deg (f)=\mp 1$ with respect
to the orientation of $S^{3}$ generate this group. The (abstract) Lie group $%
S^{3}$ defines two canonical Lie orbits $\alpha _{1},\alpha _{2}$ as
described above. Fixing some $\alpha \in \mathcal{W}(S^{3})$ arbitrarily, we
can now replace (4) by%
\begin{eqnarray}
\overline{\alpha } &:&\mathcal{W(}S^{3})/\thicksim \text{ }\longrightarrow 
\mathbb{Z} \\
&:&\alpha \longrightarrow 0  \notag
\end{eqnarray}

Now as a crucial fact, $\mathcal{H}=0$ $\Rightarrow \mathfrak{R}=0$ if $\dim
M=3$ (see Corollary 13.8, [O]). Therefore, if some $w\in \mathcal{W}$
converges, then the limit is an LLG. If $M=S^{3},$ then this LLG globalizes
to a Lie group structure since $S^{3}$ is compact and simply connected.

\begin{theorem}
Let $\alpha $ be one of the canonical orbits $\alpha _{1},\alpha _{2}\in $ $%
\mathcal{W(}S^{3})/\thicksim .$ Then all $w\in \alpha $ are convergent,
i.e., if we substitute any $w\in \alpha $ into (9) as initial condition,
then (9) converges to a limit which defines a Lie group structure on $S^{3}.$
Furthermore, all these Lie group structures are isomorphic to $S^{3}.$
\end{theorem}

Roughly, Theorem 2 states that if we start with the canonical Lie group
structure of $S^{3},$ deform it to some parallelism in the same orbit and
substitute this parallelism into HF as initial condition, then HF converges
to what we started with. Theorem 2 leaves out many questions, for instance,
if $\overline{\alpha }(\beta )=k\in \mathbb{Z}$ and if we choose $w\in \beta 
$ as initial condition and run HF, what is the behavior of the flow
depending on $k?$

Whether some $w$ integrates to a LLG is an analytical question about the
behavior of a PDE. Therefore, Theorem 2 gives an analytical characterization
of the canonical orbits $\alpha _{1},\alpha _{2}$ for $S^{3}$ and states
that what is possible homotopically, is possible also analytically. \textit{%
It is a remarkable and nontrivial fact that these canonical orbits have also
a purely homotopy theoretic characterization which is valid not only for }$%
S^{3}$ \textit{but for any homotopy }$3$\textit{-sphere. }To see this, we
should first recall some topological facts (see [KM] for the details). Let $%
N $ be a compact and oriented $4$-manifold with the (first) Pontryagin
number $p_{1}[N]$ and signature $\sigma (N).$ According to the signature
theorem we have $\sigma (N)=\frac{1}{3}p_{1}[N]$ which we write as

\begin{equation}
p_{1}[N]-3\sigma (N)=0
\end{equation}

Now let $(M,w)$ be an oriented $3$-manifold with the (positive) structure
object $w.$ Then there exists a compact and oriented $4$-manifold $N$ with
the boundary $M,$ i.e., $M=\partial N$ and the restriction of the tangent
bundle of $N$ to $\partial N=M$ is trivialized by $w.$ In this case, both
integers on the left hand side of (12) generalize to the pair $(N,\partial
N=M)$ as follows: The relative Pontryagin class $p_{1}\in H^{4}(N,\partial N;%
\mathbb{Z)}$ evaluated on the fundamental class $(N,\partial N)$ gives the
relative Pontryagin number $p_{1}[N,\partial N]$ and we also have the
relative signature $\sigma (N,\partial N)$ as the signature of the bilinear
form on $H^{2}(N,\partial N;\mathbb{R)}$ defined by cup product. However, $%
p_{1}[N,\partial N]-3\sigma (N,\partial N)=0$ does not necessarily hold in
general. The Hirzebruch defect $h(w)$ is defined in [H] by

\begin{equation}
h(w)\overset{def}{=}p_{1}[N,\partial N]-3\sigma (N,\partial N)\in \mathbb{Z}
\end{equation}

As the notation in (13) suggests, $h$ depends only on $(M,w)$ and not on $N.$
Furthermore, $h$ depends only on the orbit $[w]$ and therefore descends to a
map

\begin{equation}
h:\mathcal{W}/\thicksim \text{ }\longrightarrow \mathbb{Z}
\end{equation}

We refer to [KM] for the properties of the defect map (13). Now if $M$ has a
unique Spin structure, which is the case if $H^{1}(M,\mathbb{Z}_{2})=H^{2}(M,%
\mathbb{Z}_{2})=0,$ for instance, if $M$ is a homotopy $3$-sphere, then
there exist \textit{unique }$\beta _{1},\beta _{2}\in $ $\mathcal{W}%
/\thicksim $ \textit{satisfying}

\begin{equation}
h(\beta _{1})=2\text{ \ \ \ \ }h(\beta _{2})=-2
\end{equation}

In fact, $\beta _{1}$ and $\beta _{2}$ are the minimums of the function $%
\left\vert h(x)\right\vert ,$ $x\in \mathcal{W}/\thicksim $ (see [KM]).

\begin{definition}
Let $M$ be a homotopy $3$-sphere. Then $\beta _{1},\beta _{2}\in \mathcal{W}%
/\thicksim $ satisfying (15) are the canonical orbits.
\end{definition}

The next proposition summarizes the computation of Example 2.9 in [KM] and
states that (15) is a purely homotopy theoretic description of the canonical
Lie orbits on $S^{3}$ constructed above.

\begin{proposition}
If $M=S^{3},$ then $\beta _{1},\beta _{2}$ are the canonical Lie orbits $%
\alpha _{1},\alpha _{2}$ in Theorem 2 defined by left and right translations.
\end{proposition}

Now the following theorem generalizes Theorem 2.

\begin{theorem}
Let $M$ be a homotopy $3$-sphere and $\beta \in \mathcal{W}/\thicksim $ be a
canonical orbit. Then all $w\in \beta $ are convergent with limits
isomorphic to $S^{3}.$
\end{theorem}

Clearly PC is a consequence of Theorem 5. We hope to give the technical
details of Theorems 2, 5 in some future work.

\bigskip

\textbf{References}

\bigskip

[A] Atiyah, M. On framings of 3-manifolds, Topology, 29, (1990), 1-8

[H] Hirzebruch, F. Hilbert modular surfaces, Enseign. Math. 19 (1973),
183-281

[KM] Kirby, R., Melvin,P. Canonical framings of 3-manifolds, Turkish J. M.,
vol.23, no.1, (1999), 90-115

[O] Orta\c{c}gil, E. An Alternative Approach to Lie Groups and Geometric
Structures, OUP, 2018

\bigskip

\textbf{Acknowledgement}: I am indebted to Selman Akbulut for bringing [A]
and [KM] to my attention and explaining to me the details of some
topological constructions.

\end{document}